\documentclass[12pt]{article}

\usepackage{latexsym}
\usepackage{graphicx}
\usepackage{framed}
\usepackage{threeparttable}
\usepackage{amssymb}
\usepackage{amsmath}
\usepackage{multirow}
\usepackage{hyperref}
\usepackage{xcolor}

\begin{document}

\newcommand{\noi}{\noindent}
\newcommand{\nn}{\nonumber}
\newcommand{\bd}{\begin{displaymath}}
\newcommand{\ed}{\end{displaymath}}
\newcommand{\ds}{\displaystyle}
\newcommand{\bp}{\underline{\bf Proof}:\ }
\newcommand{\ep}{{\hfill $\Box$}\\ }

\newtheorem{1}{LEMMA}[section]
\newtheorem{2}{THEOREM}[section]
\newtheorem{3}{COROLLARY}[section]
\newtheorem{4}{PROPOSITION}[section]
\newtheorem{5}{REMARK}[section]
\newtheorem{6}{EXAMPLE}[section]
\newtheorem{7}{ALGORITHM}[section]
\newtheorem{8}{CONJECTURE}[section]

\newtheorem{10}{DEFINITION}[section]
\newtheorem{20}{OBSERVATION}[section]
\newtheorem{30}{RESULTS}[section]
\newtheorem{40}{CLAIM}[section]
\newtheorem{50}{PROBLEM}
\newcommand{\be}{\begin{equation}}
\newcommand{\ee}{\end{equation}}
\newcommand{\ba}{\begin{array}}
\newcommand{\ea}{\end{array}}
\newcommand{\bea}{\begin{eqnarray}}
\newcommand{\eea}{\end{eqnarray}}
\newcommand{\bqn}{\begin{eqnarray*}}
\newcommand{\eqn}{\end{eqnarray*}}


\newcommand{\e} { \ = \ }
\newcommand{\leqs}{ \ \leq \ }
\newcommand{\geqs}{ \ \geq \ }
\def\theequation{\thesection.\arabic{equation}}
\def\bReff#1{{\bRm
(\bRef{#1})}}
\newcommand{\eps}{\varepsilon}
\newcommand{\sgn}{\operatorname{sgn}}
\newcommand{\sign}{\operatorname{sign}}
\newcommand{\Vol}{\operatorname{Vol}}
\newcommand{\Var}{\operatorname{Var}}
\newcommand{\Cov}{\operatorname{Cov}}
\newcommand{\vol}{\operatorname{vol}}
\newcommand{\var}{\operatorname{var}}
\newcommand{\cov}{\operatorname{cov}}
\renewcommand{\Re}{\operatorname{Re}}
\renewcommand{\Im}{\operatorname{Im}}
\newcommand{\bE}{{\mathbb E}}
\newcommand{\bR}{\mathbb{R}}
\newcommand{\bN}{{\mathbb N}}
\newcommand{\bC}{\mathbb{C}}
\newcommand{\bF}{\mathbb{F}}
\newcommand{\bQ}{{\mathbb Q}}
\newcommand{\bZ}{{\mathbb Z}}
\newcommand{\cA}{{\mathcal A}}
\newcommand{\cB}{{\mathcal B}}
\newcommand{\cC}{{\mathcal C}}
\newcommand{\cD}{{\mathcal D}}
\newcommand{\cE}{{\mathcal E}}
\newcommand{\cF}{{\mathcal F}}
\newcommand{\cG}{{\mathcal G}}
\newcommand{\cH}{{\mathcal H}}
\newcommand{\cI}{{\mathcal I}}
\newcommand{\cJ}{{\mathcal J}}
\newcommand{\cK}{{\mathcal K}}
\newcommand{\cL}{{\mathcal L}}
\newcommand{\cM}{{\mathcal M}}
\newcommand{\cN}{{\mathcal N}}
\newcommand{\cO}{{\mathcal O}}
\newcommand{\cP}{{\mathcal P}}
\newcommand{\cQ}{{\mathcal Q}}
\newcommand{\cR}{{\mathcal R}}
\newcommand{\cS}{{\mathcal S}}
\newcommand{\cT}{{\mathcal T}}
\newcommand{\cU}{{\mathcal U}}
\newcommand{\cV}{{\mathcal V}}
\newcommand{\cW}{{\mathcal W}}
\newcommand{\cX}{{\mathcal X}}
\newcommand{\cY}{{\mathcal Y}}
\newcommand{\cZ}{{\mathcal Z}}
\newcommand{\bx}{{\mathbf x}}
\newcommand{\by}{{\mathbf y}}
\newcommand{\bz}{{\mathbf z}}
\newcommand{\bba}{{\mathbf a}}
\newcommand{\bbb}{{\mathbf b}}
\newcommand{\bbc}{{\mathbf c}}


\title{Marginal Stationary Distributions and Convergence Rates of Higher Order Markov Chains}
\author{
Lixing Han\thanks{Department of Mathematics, University of Michigan-Flint, Flint, MI 48502, USA. Email: \texttt{lxhan@umich.edu}} 
}
\date{}
\maketitle

\begin{abstract}  
 For a regular higher order Markov chain, its reduced first order chain may be reducible. When the reduced first order chain is reducible, it does not have a unique stationary distribution. However, it has been shown that the probability distribution of the current state of a regular higher order chain converges to a unique limiting probability distribution. In this paper, we interpret this limiting distribution as a marginal stationary distribution and prove some properties of marginal stationary distributions. We then establish some convergence rate results for such type of convergence. Finally, we introduce two types of marginal mixing times, which extend the notions of mixing times for first order Markov chains to higher order chains. 
\end{abstract}

\ \\
{\bf Key words.}  Higher order Markov chain,   marginal stationary distribution, marginal convergence rate, marginal mixing time

\ \\
{\bf AMS subject classification (2010). 60J10, 15A69}

\section{Introduction}
\label{Intro}
\setcounter{equation}{0}

 Rates of convergence and mixing times are a central object in the theory of first order Markov chains with important applications in several areas, including statistical physics and theoretical computer science. In this paper, we consider these topics for higher order Markov chains.
 In recent years, higher order Markov chains  have been a subject of much research (\cite{BozHaj16, ChaZha13, CulPeaZha17, DinNgWei18, FasTud20, Gei17, GleLimYu15, HanWanXu22, HanXu24a, HanXu24b, HanXu26, HuQi14, HuaQi20, LiZha16, LiNg14, WuChu17}). 
 
Let $m \ge 2$ and $n \ge 1$ be integers. A stochastic chain   $\{X_t\}_{t=1}^{\infty}$ with the state space $S=\{1,2, \ldots, n\}$ is said to be an $(m-1)$th order Markov chain if it satisfies the following property
\be
\label{tran-1}
\Pr(X_t=i_1| X_{t-1}=i_2, \ldots, X_{1}=i_t) =\Pr(X_t=i_1| X_{t-1}=i_2, \ldots, X_{t-m+1}=i_{m}),
\ee
for $i_1, i_2, \ldots, i_{t} \in S$ and  $t \ge m$.    In this paper, we will focus on homogeneous chains, in which the transition probabilities (\ref{tran-1}) are independent of the time $t$. Denote
\be
p_{i_1 i_2 \ldots i_{m}} = \Pr(X_t=i_1| X_{t-1}=i_2, \ldots, X_{t-m+1}=i_{m}).
\ee 
Then ${\cal P}=[p_{i_1 i_2 \ldots i_{m}}]$ is called the transition tensor of the $(m-1)$th order Markov chain, which is an $m$th order $n$-dimensional stochastic tensor satisfying $p_{i_1 i_2 \ldots i_{m}} \ge 0$ for  $i_1, i_2, \ldots, i_{m} \in S$, and
$$
\sum_{i_1 \in S}  p_{i_1 i_2, \ldots,i_{m}}=1, \ \  \forall i_2, \ldots, i_{m} \in S.
$$
Denote the probability distribution of the chain at time $t$ by $\xi^{(t)} \in \mathbb{R}^n$, whose $i$th entry is  
\be
\label{prob-distr}
\xi^{(t)}_{i} = \Pr(X_t=i).
\ee

The $k$-step, with $k \ge 1$, transition probability from states $(i_2,\ldots,i_m)$ to state $i_1$ is given by 
$$
p^{(k)}_{i_1i_2\ldots i_m}=\Pr(X_{k+m-1}=i_1 | X_{m-1}=i_2, \ldots, X_{1}=i_m).
$$
 Set ${\cal P}^{(k)}=[p^{(k)}_{i_1\ldots i_m}]$.    ${\cal P}^{(k)}$ is also an $m$th order $n$-dimensional stochastic tensor and is called the $k$-step transition tensor. Clearly, ${\cal P}^{(1)} = {\cal P}$ and
 $$
p^{(k+1)}_{i_1i_2\ldots i_m}=\sum_{j \in S} p_{ji_2\ldots i_m}p^{(k)}_{i_1ji_2\ldots i_{m-1}}, ~k=1, 2, \ldots.
$$

A central issue in the study of Markov chains is to investigate their long term behavior, including the existence and uniqueness of stationary distributions, and the convergence of distributions $\xi^{(t)}$ to the stationary distribution and the rates of convergence.  

When $m=2$, the chain $\{X_t\}_{t=1}^{\infty}$  is a standard first order Markov chain with the transition matrix $P=(p_{ij})$. It is well known that (see, for example, \cite{Ios07, KemSne60, LevPer17}),   if the first order chain is regular (that is, irreducible and aperiodic), then the  distribution $\xi^{(t)}$
converges to the unique stationary distribution $\pi \in \mathbb{R}^n$ from any initial distribution $\xi^{(1)}$, where $\pi$ satisfies  $\pi_i \ge 0$ for $i \in S$, $\sum_{i \in S} \pi_i=1$, and 
$$
\pi_i = \sum_{j \in S} P_{ij} \pi_j,   \ \   \forall i \in S.
$$

When $m \geq 3$, an $(m-1)$th order chain is called a higher order Markov chain. The stationary distributions of higher order  Markov chains are  defined via its reduced first order chain. Let $Y_t = (X_{t+m-2},  \ldots, X_{t})$. Then  the higher order chain 
$\{X_t\}_{t=1}^{\infty}$  can be reduced to a first order chain $\{Y_t\}_{t=1}^{\infty}$ on the product space $S^{m-1}$ (see, for example, Hunter \cite{Hun83}). The probability distribution of $Y_t$ can be written as an $(m-1)$th-order $n$-dimensional 
tensor ${\cal Z}^{(t)} = [z_{i_1 i_2 \ldots i_{m-1}}^{(t)}]$, where
$$
z^{(t)}_{i_1 i_2 \ldots i_{m-1}} = \Pr ( X_{t+m-2}=i_1, X_{t+m-3}=i_2, \ldots, X_{t}=i_{m-1}), \ \ \ i_1, \ldots, i_{m-1} \in S. 
$$
Note that  ${\cal Z}^{(t)}$ is a joint distribution of $X_{t+m-2},  \ldots, X_{t}$. The stationary distributions of  a higher order chain are defined as follows.

\begin{10}  {\rm
 A stationary distribution of an  $(m-1)$th order Markov chain 
$\{X_t\}_{t=1}^{\infty}$ on $S$ with transition tensor $\cal P$ is an $(m-1)$th-order $n$-dimensional tensor $\mathcal{W}=[w_{i_1 i_2 \ldots i_{m-1}} ] $ satisfying
\be
\label{stationary-distribution}
w_{i_1 i_2 \ldots i_{m-1}} = \sum_{i_{m} \in S} p_{i_1 i_2 \ldots i_{m}} w_{i_2 i_3 \ldots i_{m}}, \ \ \sum_{i_1, i_2,  \ldots, i_{m-1}  \in S} w_{i_1 i_2 \ldots i_{m-1}} =1, \ \ w_{i_1 i_2 \ldots i_{m-1}} \ge 0,
\ee
for all  $ i_1, \ldots, i_{m-1} \in S$.
}
\end{10}

Clearly, $\mathcal{W}$ is a stationary distribution of $\{X_t\}_{t=1}^{\infty}$ if and only if  $\mathcal{W}$ is a stationary distribution of $\{Y_t\}_{t=1}^{\infty}$.  If ${\cal Z}^{(1)}={\cal W}$, then
$$
{\cal Z}^{(t)}={\cal W}, 
$$
for $t=1,2,\ldots$. As $Y_t = (X_{t+m-2},  \ldots, X_{t})$, we will call $\mathcal{W}$ a joint stationary distribution in this paper. 

Using a proper indexing (see Section 2), we can formulate an $n^{m-1} \times n^{m-1}$ transition matrix $Q$ for the reduced first order chain $\{Y_t\}_{t=1}^{\infty}$, and accordingly vectorize ${\cal Z}^{(t)}$ and $\mathcal{W}$ into probability vectors $\zeta^{(t)}$ and $\omega$ in $\mathbb{R}^{n^{m-1}}$,  respectively. With these notations, Condition (\ref{stationary-distribution}) is equivalent to
\be
\label{Qw}
\omega=Q \omega, \ \ \sum_{i=1}^{n^{m-1}} \omega_i=1, \ \ \omega_i \ge0 \ ( i=1,2, \ldots, n^{m-1}).  
 \ee
The questions on existence and uniqueness of joint stationary distributions of a higher order chain $\{X_t\}_{t=1}^{\infty}$, the convergence of the joint distributions 
$\{{\cal Z}^{(t)}\}$ to the joint stationary distribution, and the convergent rate can therefore be answered using the known results for first order chains, provided that the reduced first order chain $\{Y_t\}_{t=1}^{\infty}$ is ergodic or regular. 

However, it has been discovered that problems regarding a higher order chain  $\{X_t\}_{t=1}^{\infty}$ cannot always be fully addressed by resorting to its reduced first order chain $\{Y_t\}_{t=1}^{\infty}$. For example, $\{Y_t\}_{t=1}^{\infty}$ may be non-ergodic even though$\{X_t\}_{t=1}^{\infty}$ is ergodic \cite{HanXu24a}.

Regularity of first order Markov chains has been extended to higher order chains (\cite{HanXu26, Vla85}).   In  \cite{HanXu26}, it is shown that the reduced first order chain of a regular higher order chain  may be non-regular (in fact, it can  be reducible).  
When a higher order chain $\{X_t\}_{t=1}^{\infty}$  is regular but its reduced first order chain $\{Y_t\}_{t=1}^{\infty}$ is reducible,  the joint stationary distributions are not unique and the distribution of $Y_t$ may not be convergent.  In this case,  the long term behavior of the higher order chain $\{X_t\}_{t=1}^{\infty}$ cannot be obtained via the convergence of the distribution of $Y_t$ . However,  it is proved in \cite{HanXu26} that if a higher order chain $\{X_t\}_{t=1}^{\infty}$ is regular, then  $\xi^{(t)}$, the probability distribution of the chain at time $t$,  converges to a unique limiting probability distribution $\pi \in \mathbb{R}^n$ from any initial distribution of $(X_{m-1}, \ldots, X_{1})$, even if its reduced first order chain is reducible. This type of convergence can be used to predict the long term behavior of a regular higher order chain. To better understand the long term behavior, the following questions arise naturally:  

\begin{itemize}
 \item What does the limiting distribution $\pi$ mean? 
 \item How is $\pi$ related to the joint stationary distributions of the chain $\{X_t\}_{t=1}^{\infty}$? 
 \item How fast does the sequence $\{\xi^{(t)}\}$ converge to $\pi$?
 \end{itemize}

In this paper, we answer these questions. We provide an interpretation of the limiting distribution, which is a marginal stationary distribution. We also analyze the rate of convergence of  $\{\xi^{(t)}\}$.

This paper is organized as follows. In Section 2, we summarize some known results about regular or ergodic higher order chains which will be used in later sections. In Section 3, we introduce marginal stationary distributions and prove their properties. In Section 4, we establish some convergence rate results for  convergence to the marginal stationary distribution. Finally in Section 5,  we introduce two notions of marginal mixing times.

\section{Preliminaries}
\label{Prelim}
\setcounter{equation}{0}

In this section, we review some concepts and results about higher order Markov chains that will be used in the next few sections. We assume $m \ge 3$. 

\begin{10} {\rm
An $(m-1)$the order Markov chain with transition tensor $\cal P$ is siad to be {\it ergodic},\footnote{We note that in the first order case, an ergodic chain is also called an irreducible chain and a regular chain is equivalent to an irreducible and aperiodic chain. However, for higher order chains, irreducibility and ergodicity differ (see, for example, \cite{HanWanXu22}).} if given any $i_1, i_2, \ldots, i_m \in S$, there exists $k \ge 1$, which may depend on $i_1, i_2, \ldots, i_m$, such that $p^{(k)}_{i_1i_2\ldots i_m}>0$.  

It is called  {\it regular} if there exists some fixed $h \ge 1$ such that $p^{(h)}_{i_1i_2\ldots i_m} > 0$ for every $i_1, i_2, \ldots, i_m \in S$.
 }
\end{10}

 Denote the first passage time random variable from states $(i_2,\ldots,i_m)$ to state $i_1$ by
 \be
 \label{first-passage}
 T_{ i_1 i_2 \ldots i_m } = \min \{ j \geq 1: X_{m-1+j}=i_1| X_{m-1}=i_2, \ldots, X_{1}=i_{m} \}. 
 \ee
   The mean first passage time is defined as
 $$
 \mu_{ i_1 i_2 \ldots i_m } = E (T_{ i_1 i_2 \ldots i_m }  ),
 $$
 wehre $E(\cdot)$ denotes the expected value. Regarding mean first passage times, we have the following result.
 
 \begin{2} {\rm \cite{HanXu24a}}
 \label{MFPT}
  For an $(m-1)$th order ergodic Markov chain,  the mean first passage times $\mu_{ i_1 i_2 \ldots i_m }$ are finite for each $i_1, i_2, \ldots, i_m \in S$ and satisfy
 \be
 \label{sys}
 \mu_{ i_1 i_2 \ldots i_m} =1 + \sum_{j \ne i_1}  p_{j i_2 \ldots i_m} \mu_{ i_1 j i_2 \ldots i_{m-1}}.
 \ee
 Moreover, the linear system (\ref{sys}) is nonsingular  if and only if the chain is ergodic.
 \end{2}

\begin{2}
\label{conv} {\rm \cite{HanXu26}}
Suppose that  $\{X_t \}_{t=1}^{\infty}$ is  an $(m-1)$th order Markov chain on the state space $S$ with a transition tensor $\cal P$. If the chain is regular, then there exists 
$\pi \in \mathbb R^n$ satisfying $\pi_i > 0$ for each $i \in S$ and $\sum_{i=1}^n \pi_i = 1$ such that 
\be
\label{lim_for}
\lim_{k \to \infty} p^{(k)}_{ii_2\ldots i_m}=\pi_i, \ \ \ \  \forall i_2, \ldots, i_m \in S.
\ee 
Moreover, 
\be
\label{prob}
\lim_{t \to \infty} \xi_{i}^{(t)} = \lim_{t \rightarrow \infty}\Pr(X_t=i)=\pi_i, \ \ \ \  \forall i \in S,
\ee
regardless of the initial probability distributions of $(X_{m-1}, \ldots, X_1)$. 
\end{2}

We now review how to reduce an $(m-1)$th order chain $\{X_t\}_{t=1}^{\infty}$ to a first order chain $\{Y_t\}_{t=1}^{\infty}$, where  $Y_t = (X_{t+m-2},  \ldots, X_{t})$. More detailed derivations can be found in \cite{HanXu26}.  Denote the product space $S^{m-1}=\{i_1i_2\ldots i_{m-1} : i_1, i_2, \ldots, i_{m-1} \in S\}$ as a collection of multi-indices of length $m-1$. These multi-indices are arranged via linear indexing \cite{MarShaLar13}. As an example, when $m=4$ and $n=2$,  we have
$S^3=\{111, 211, 121, 221, 112, 212, 122, 222\}.$
With this notation, when $X_{t+m-2}=i_1,  \ldots, X_{t}=i_{m-1},$ 
we have $Y_t=i_1i_2\ldots i_{m-1}$.  Let $N=n^{m-1}$. The transition matrix $Q \in \mathbb R^{N \times N}$ of $\{Y_t\}_{t=1}^{\infty}$ is given by 
$$
q_{i_1i_2\ldots i_{m-1}, j_2j_3\ldots j_m}=\Pr(Y_{t+1}=i_1i_2\ldots i_{m-1} \ | \ Y_t=j_2j_3\ldots j_m).
$$
Note that  $q_{i_1i_2\ldots i_{m-1}, j_2j_3\ldots j_m}=p_{i_1i_2\ldots i_{m-1}j_m}$ when $i_\ell=j_\ell$ for all $\ell=2, 3, \ldots, m-1$; otherwise, $q_{i_1i_2\ldots i_{m-1}, j_2j_3\ldots j_m}=0$.

Let $P \in \mathbb R^{n \times N}$ be the mode-$1$ matricization \cite{KolBad09} of the transition tensor ${\cal P}$. That is, $P$ consists of the frontal slices of ${\cal P}$ which are arranged side by side and ordered according to linear indexing. 
 Define the $n \times N$ matrix
\be
\label{P0}
P^{(0)}=[\underbrace{I_n ~I_n ~\ldots ~I_n}_{n^{m-2}}],
\ee
where $I_n \in \mathbb R^{n \times n}$ is the identity matrix. Then an alternative way to get the transition matrix  $Q$ from ${\cal P}$ is   by 
$$
Q=P^{(0)} \ast P,
$$
where $\ast$ stands for the columnwise Khatri-Rao product \cite{SmiBroGel04}. 

On the other hand, the $k$-step transition tensor ${\cal P}^{(k)}$ can be obtained from $Q^k$, which is the $k$th power of $Q$ (i.e.,  the $k$-step transition matrix of $\{Y_t\}_{t=1}^{\infty}$). Let $P^{(k)} \in \mathbb R^{n \times N}$ be the mode-$1$ matricization of the $k$-step transition tensor ${\cal P}^{(k)}$. Then we have
\be
\label{PkQk}
P^{(k)}=P^{(0)}Q^k,
\ee
for $k=1,2, \ldots$.

From (\ref{PkQk}), if $Q$ is regular, then ${\cal P}$ is also regular. However, the converse is not true. \cite{HanXu26} gives an example in which  
${\cal P}$ is also regular, but $Q$ is reducible.

\begin{2} {\rm \cite{HanXu26}}
\label{conn}
Suppose that $\{X_t\}_{t=1}^{\infty}$ is an $(m-1)$th order regular Markov chain on state space $S$. Let $Q$ be the transition matrix of its reduced first order chain $\{Y_t\}_{t=1}^{\infty}$. If $\omega \in \mathbb R^N$ is a normalized nonnegative right eigenvector of $Q$ associated with dominant eigenvalue $1$, i.e. $Q\omega=\omega$, $\omega_i \ge 0$ for $i=1,2, \ldots, N$, and $\sum_{i=1}^{N} \omega_i=1$. Then,
\be
\label{pz}
\pi = P^{(0)}\omega
\ee 
is a probability distribution in $\mathbb{R}^n$, which is a constant vector for any choice of $\omega$. Moreover, this $\pi$ is the unique limiting probability distribution satisfying (\ref{prob}). 
\end{2}

\section{Marginal stationary distributions}
\label{Stationary-Distr}
\setcounter{equation}{0}

We start with the existence and uniqueness of joint stationary distributions of a higher order chain. Suppose that $m \ge 3$. 

\begin{4}
\label{existence-uniqueness}
  An $(m-1)$th order Markov chain $\{X_t\}_{t=1}^{\infty}$ on state space $S$ with a transition tensor  $\cal P$ always has a joint stationary distribution.  If the reduced first order chain $\{Y_t\}_{t=1}^{\infty}$ with transition matrix $Q$ is ergodic (i.e., $Q$ is irreducible), then  $\{X_t\}_{t=1}^{\infty}$ has a unique joint stationary distribution. 
\end{4}
\bp  $\{X_t\}_{t=1}^{\infty}$ has a joint stationary distribution if and only if (\ref{Qw}) has a solution. Since $Q$ is a stochastic matrix, by the Perron-Frobenus theorem, (\ref{Qw}) always has a solution. Thus, $\{X_t\}_{t=1}^{\infty}$ always has a joint stationary distribution. If $Q$ is irreducible, then  (\ref{Qw}) has a unique solution, implying $\{X_t\}_{t=1}^{\infty}$ has a unique joint stationary distribution.
\ep

\begin{5} {\rm
When $m \ge 3$, under what conditions on the transition tensor $\cal P$ of an $(m-1)$th order chain $\{X_t\}_{t=1}^{\infty}$ that its reduced first order chain  $\{Y_t\}_{t=1}^{\infty}$ is ergodic is a challenging question and currently has only partial answers.  
In \cite{HanXu26}, it is shown that when $\mathcal{P}$ is positive,  $Q$ is regular (therefore, ergodic). Geiger \cite{Gei17} gives a sufficient condition for the existence of a unique joint stationary distribution of a higher order chain. However, the condition in \cite{Gei17} is sophisticated and it is unknown whether the distribution of $Y_t$ converges to the joint stationary distribution under this condition.   
}
\end{5}

\begin{10} 
\label{marginal-def}{\rm 
If $\mathcal{W}$ is a joint stationary distribution of an $(m-1)$th order Markov chain  $\{X_t\}_{t=1}^{\infty}$ on state space $S$, then $\pi \in \mathbb{R}^n$
defined by
\be
\label{marginal-distr}
\pi_i = \sum_{i_2,\ldots, i_{m-1} \in S} w_{i i_2 \ldots i_{m-1}}
\ee
is called a marginal stationary distribution of the chain.
}
\end{10}

By Proposition \ref{existence-uniqueness} and Definition \ref{marginal-def}, we have the following existence result.

\begin{3}{\rm 
 An $(m-1)$th order Markov chain always has a marginal stationary distribution. If its reduced first order chain is ergodic, then it has a unique marginal  stationary distribution. 
}
\end{3}

The following proposition shows that marginal stationary distributions satisfies a time invariance property, that is, the marginal stationary distribution are the same at times $t+m-2, t+m-3, \ldots, t$.  

\begin{4}
\label{property1}
Suppose that $\mathcal{W}$ is a joint stationary distribution of an $(m-1)$th order Markov chain  $\{X_t\}_{t=1}^{\infty}$ on state space $S$ with a transition tensor $\mathcal{P}$. If $\pi$ is the marginal stationary distribution defined by $\mathcal{W}$, then we have
$$
\pi_i=\sum_{i_2, i_3, \ldots, i_{m-1} \in S} w_{i i_2 \ldots i_{m-1}} = \sum_{i_1, i_3, \ldots, i_{m-1} \in S} w_{i_1 i i_3 \ldots i_{m-1}} = \cdots = \sum_{i_1,  \ldots, i_{m-2} \in S} w_{i_1  \ldots i_{m-2} i}.
$$
\end{4}
\bp By (\ref{stationary-distribution}), we have
\[
\begin{array}{rcl}
\displaystyle{\sum_{i_1, i_{3},  \ldots, i_{m-1} \in S} w_{i_1 i i_3  \ldots i_{m-1}}}
&=& \displaystyle{ \sum_{i_1, i_{3},  \ldots, i_{m-1} \in S} \left ( \sum_{i_{m} \in S}  p_{i_1 i i_3 \ldots i_{m}} w_{i i_3 \ldots i_{m}} \right )} \\
& = & \displaystyle{ \sum_{i_{3},  \ldots,  i_{m} \in S} \left ( \sum_{i_1 \in S} p_{i_1 i \ldots i_{m}} \right ) w_{i i_3 \ldots i_{m}}} \\
& = & \displaystyle{ \sum_{i_{3},  \ldots, i_{m} \in S } w_{i i_3 \ldots i_{m}}} \\
& = & \pi_{i}.
\end{array}
\]
The other equalities can be proved similarly. 
\ep

The following proposition shows that if a higher order chain  $\{X_t\}_{t=1}^{\infty}$ starts with an joint stationary distribution $\mathcal{W}$, then the distribution of $\{Y_t\}_{t=1}^{\infty}$ is invariant for any $t$. Moreover, the marginal distribution of $X_{t+m-2}$ in $Y_t$ is also time invariant and equals to the marginal distribution defined by $\mathcal{W}$.   Its proof is straightforward and omitted.  
\begin{4}
\label{property2}
Suppose that $\mathcal{W}$ is a joint stationary distribution of an $(m-1)$th order Markov chain  $\{X_t\}_{t=1}^{\infty}$ on state space $S$ with transition tensor $\mathcal{P}$ and $\pi$ is the  marginal stationary distribution defined by $\mathcal{W}$. If the chain starts with $Y_1=(X_{m-1}, \ldots, X_1)$ that has the joint distribution $\mathcal{Z}^{(1)}=\mathcal{W}$, then the joint distribution ${\cal Z}^{(t)}$ of $Y_t=(X_{t+m-2}, \ldots, X_t)$ satisfies
$$
\mathcal{Z}^{(t)} = \mathcal{W}
$$
for all $t \ge 1$. Furthermore, let $\beta^{(t)} \in \mathbb{R}^n$ be the marginal distribution of $X_{t+m-2}$ in $Y_t$, that is, 
$$
\beta^{(t)}_i=\sum_{i_2, \ldots, i_{m-1} \in S} z^{(t)}_{i i_2 \ldots i_{m-1}}, \ \ i \in S.
$$
Then we have
$$
\beta^{(t)} = \pi, 
$$
for all $t \ge 1$.
\end{4}

The next theorem shows that a regular higher order Markov chain must have a unique marginal stationary distribution.

\begin{2}
\label{uniqueness}
Suppose that $\{X_t\}_{t=1}^{\infty}$ is an $(m-1)$th order Markov chain   on state space $S$ with a transition tensor $\mathcal{P}$. If the chain is regular, then it has a unique marginal stationary distribution $\pi$. Moreover, 
$$
\lim_{t \to \infty} \xi^{(t)}_{i} =  \lim_{t \to \infty} \Pr(X_t=i) =  \pi_i, \ \ \forall i \in S. 
$$
Thus, $\pi$ is the limiting probability distribution of $X_t$. 
\end{2}
\bp
By the equivalence of (\ref{stationary-distribution}) and  (\ref{Qw}), all marginal stationary distributions $\pi$ can be obtained by
$$
\pi= P^{(0)} \omega, 
$$
where $P^{(0)}$ is defined in (\ref{P0}) and $\omega$ is a solution of  (\ref{Qw}). By Theorem \ref{conn}, if  $\{X_t\}_{t=1}^{\infty}$ is regular, there is a unique marginal stationary distribution $\pi$,  and  $\lim_{t \to \infty} \xi^{(t)}_{i}  =  \pi_i, \ \ \forall i \in S$.
\ep

We are now ready to an interpretation of the limiting distribution obtained in Theorem \ref{conv}. For a regular higher order Markov chain, its  
limiting distribution is its unique marginal stationary distribution.

\begin{5} {\rm We comment that a regular higher order Markov chain may have more than one joint stationary distributions. Consider the second order chain with  the following transition tensor (\cite{Gei17}):
$${\cal P}(:,:,1)=\left[\ba{cccc}
1/2 & 0 & 0 & 0\\
1/2 & 0 & 1 & 0\\
0 & 1 & 0 & 1\\
0 & 0 & 0 & 0
\ea\right], ~{\cal P}(:,:,2)=\left[\ba{cccc}
0 & 0 & 1/2 & 1\\
0 & 1/2 & 0 & 0\\
1/2 & 1/2 & 0 & 0\\
1/2 & 0 & 1/2 & 0
\ea\right],$$
$${\cal P}(:,:,3)=\left[\ba{cccc}
0 & 1 & 0 & 1\\
1 & 0 & 1/2 & 0\\
0 & 0 & 1/2 & 0\\
0 & 0 & 0 & 0
\ea\right], ~{\cal P}(:,:,4)=\left[\ba{cccc}
0 & 0 & 0 & 0\\
1 & 1 & 1 & 0\\
0 & 0 & 0 & 1/2\\
0 & 0 & 0 & 1/2
\ea\right].$$
It is shown in \cite {HanXu26} that this chain is regular, but its reduced first order chain is reducible. This chain has a unique marginal stationary distribution  
$$\pi=[2/7, 2/7, 2/7, 1/7]^T.
$$
 However, it has two different joint stationary distributions:
$$
W_1=\left[\ba{cccc}
0 & 2/7 & 0 & 0\\
0 & 0 & 2/7 & 0\\
1/7 & 0 & 0 & 1/7\\
1/7 & 0 & 0 & 0
\ea\right],  \ \ \ \ W_2=\left[\ba{cccc}
0 & 0 & 1/7 & 1/7\\
2/7 & 0 & 0 & 0\\
0 & 2/7 & 0 & 0\\
0 & 0 & 1/7 & 0
\ea\right].
$$ 
In fact, any convex combination of $W_1$ and $W_2$ is also a joint stationary distribution.  
}
\end{5}

We end this section by giving a discussion on joint and marginal stationary distributions.  Consider an $(m-1)$th order Markov chain. The joint stationary distribution gives the long-run probability of finding the system in a  sequence of states over $m-1$ consecutive steps. On the other hand, the marginal stationary distribution  focuses on a single point in time. The marginal distribution $\pi_i$ tells the long-run probability that the system is in state $i$ at a single point in time. The  marginal stationary distribution  can play a significant role in predicting the long term behavior of a higher order chain when it is regular but its reduced first order chain is non-regular.

\section{Convergence rates}
\label{Rates}
\setcounter{equation}{0}

We first establish a geometric rate of convergence result for regular higher order Markov chains. 

\begin{2}
\label{rate1}
Suppose that $\{X_t\}_{t=1}^{\infty}$ is an $(m-1)$th order regular Markov chain  on state space $S$ with a transition tensor $\mathcal{P}$.  Let $\pi$ be its unique marginal stationary distribution. Then there exist constants $C >0$ and $ 0<D<1$ such that
\be
\label{geometric1}
\max_{i_2,\ldots, i_m} | p_{i i_2 \ldots i_m}^{(k)} - \pi_i | \le C D^k,
\ee
for all $i \in S$ and $k \ge 1$. 
\end{2}

\bp
To prove the theorem, we adopt the technique used in the proof of Theorem 3.3 in \cite{HanXu26}. 

We first consider the case when ${\cal P}$ is positive. Let
$$
\epsilon = \min_{i_1,\ldots,i_m} p_{i_1 \ldots i_m}. 
$$

Note that $\epsilon > 1/2$ occurs only when $n=1$, and  $\epsilon = 1/2$ occurs only when $n=2$ and  $p_{i_1 \ldots i_m} \equiv 1/2$ for all 
$i_1,\ldots,i_m \in S$. In both cases, (\ref{geometric1}) trivially holds.  WLOG, assume that $0<\epsilon < 1/2$.

Fix $i \in S$.  Let $x$ be the $i$-th standard basis vector $e_i \in \mathbb{R}^n$, i.e., $x_i = 1$ and $x_j = 0$ for all $j \neq i$. According to Equation (3.5) in \cite{HanXu26}, define 
$$y^{(k)}_{i_2 \ldots i_m} = \sum_{j=1}^n x_j p^{(k)}_{j i_2 \ldots i_m}.
$$
 For the  chosen $x$, this simplifies to $y^{(k)}_{i_2 \ldots i_m}= p^{(k)}_{i i_2 \ldots  i_m}$.
Let $U_0=\max_i \{x_i\}=1$, $L_0=\min_i \{x_i\} = 0$. Set 
\be
\label{Uk-Lk}
U_k = \max_{i_2 \ldots i_m} y^{(k)}_{i_2 \ldots i_m}= \max_{i_2...i_m} p^{(k)}_{i i_2 \ldots i_m}, \ \ 
L_k = \min_{i_2, \ldots, i_m} y^{(k)}_{i_2 \ldots i_m} = \min_{i_2, \ldots, i_m} p^{(k)}_{i i_2 \ldots i_m},
\ee
for $ k \ge 1$, where $p^{(k)}_{i i_2 \ldots i_m}$ is the $(i, i_2, \ldots, i_m)$ entry of ${\cal P}^{(k)}$.

By the proof of Theorem 3.3 in \cite{HanXu26}, we have
$$
U_1-L_1 \le (1-2 \epsilon ) (U_0-L_0). 
$$
Moreover, the sequence $\{U_k\}$ is decreasing and $\{L_k\}$ is increasing,  that is,
$$
U_{k+1} \le U_k, \ \ \  L_{k+1} \ge L_k.
$$ 
Thus, the sequence of differences $\{U_k - L_k\}$ is decreasing. The proof of Theorem 3.3 in \cite{HanXu26} establishes the following inequalities for $k \ge m - 1$
$$
U_{k+1} \le U_{k-m+2} - \epsilon^{m-1}(U_{k-m+2} - L_{k-m+2}), 
$$
and
$$
 L_{k+1} \ge L_{k-m+2} + \epsilon^{m-1}(U_{k-m+2} - L_{k-m+2}).
$$
Subtracting the second inequality from the first yields
\be
\label{recursive}
U_{k+1} - L_{k+1} \le (1 - 2\epsilon^{m-1})(U_{k-m+2} - L_{k-m+2}).
\ee

Let $\alpha = 1 - 2\epsilon^{m-1}$. Since $0 < \epsilon < 1/2$, we have $0 < \alpha < 1$. Applying (\ref{recursive}) iteratively for $q$ steps, where each step jumps by $m-1$ indices, we obtain
$$
U_{(m-1)q} - L_{(m-1)q} \le \alpha^q (U_0 - L_0).
$$
For any integer $k \ge 1$, we can write $k = q(m-1) + r$, where $0 \le r < m-1$. 
As $\{U_k - L_k\}$ is a decreasing sequence, we have
\be
\label{upperbd}
U_k - L_k \le U_{q(m-1)} - L_{q(m-1)} \le \alpha^q (U_0 - L_0).
\ee
Note that $q = \frac{k-r}{m-1} > \frac{k}{m-1} - 1$. Substituting this into (\ref{upperbd}) gives
$$
U_k - L_k \le \alpha^{\frac{k}{m-1} - 1} (U_0 - L_0) = \alpha^{-1} \cdot (\alpha^{\frac{1}{m-1}})^k \cdot (U_0 - L_0).
$$

Let 
$$
D = \alpha^{\frac{1}{m-1}} = (1 - 2\epsilon^{m-1})^{\frac{1}{m-1}}, \ \ C= \alpha^{-1} (U_0 - L_0)= \frac{1}{1 - 2\epsilon^{m-1}}.
$$
 Since $0 < \alpha < 1$, it follows that $0 < D < 1$ and $C>0$.  This leads to
 $$
 U_k - L_k \le C \cdot D^k.
 $$
 Note that $C$ and $D$ are independent of $i$.

By Theorem \ref{conv}, $\lim_{k\to\infty} p^{(k)}_{i i_2...i_m} = \pi_i$. Since $L_k \le p^{(k)}_{i i_2...i_m} \le U_k$ for all $k \ge 1$, we have $L_k \le \pi_i \le U_k$ for all $k \ge 1$. Therefore, we have
$$
|p^{(k)}_{i i_2...i_m} - \pi_i| \le U_k - L_k \le C \cdot D^k,
$$
for all $i, i_2, \ldots, i_m \in S$ and $k \ge 1$.  This implies that (\ref{geometric1}) holds.  \\

If the transition tensor $P$ is regular but not strictly positive, by definition there exists some integer $h > 1$ such that the $h$-step transition probability tensor is strictly positive, i.e., ${\cal P}^{(h)} > 0$. We can find the constants $C$ and $D$ by adapting the above proof to operate over blocks of $h$ steps.

Let
$$
\epsilon = \min_{i_1, i_2 \ldots i_m} p^{(h)}_{i_1 i_2 \ldots i_m}.
$$
Aagin, we can assume $0 < \epsilon < 1/2$. Let Let $U_0=1$, $L_0 = 0$, and $U_k$ and $L_k$ be defined as in (\ref{Uk-Lk}).

 As ${\cal P}^{(h)}$ acts as a strictly positive transition tensor, the contraction of the bounds $U_k$ and $L_k$ happens over $m-1$ applications of 
 ${\cal P}^{(h)}$. Using a similar argument, we have
 $$
 U_{k+h(m-1)} - L_{k+h(m-1)} \le (1 - 2\epsilon^{m-1})(U_k - L_k).
 $$
 
 Let $\alpha = 1 - 2\epsilon^{m-1}$. Then $0 < \alpha < 1$. Write
$$k = q \cdot h(m-1) + r, 
$$
where  $0 \le r < h(m-1)$. Because $\{U_k - L_k\}$ is a decreasing sequence, we have
$$
U_k - L_k \le U_{q \cdot h(m-1)} - L_{q \cdot h(m-1)} \le \alpha^q (U_0 - L_0).
$$
Since $q = \frac{k - r}{h(m-1)} > \frac{k}{h(m-1)} - 1$, we obtain
$$
U_k - L_k \le \alpha^{\frac{k}{h(m-1)} - 1} (U_0 - L_0) = \alpha^{-1} \cdot \left(\alpha^{\frac{1}{h(m-1)}}\right)^k \cdot (U_0 - L_0).
$$

Let
$$
D = \alpha^{\frac{1}{h(m-1)}} = (1 - 2\epsilon^{m-1})^{\frac{1}{h(m-1)}}, \ \ C = \alpha^{-1} (U_0 - L_0) = \frac{1}{1 - 2\epsilon^{m-1}}
$$
Then $0<D<1$ and $C>0$, and  (\ref{geometric1}) holds. 
\ep

\begin{3}
\label{rate2}
Suppose that $\{X_t\}_{t=1}^{\infty}$ is an $(m-1)$th order regular Markov chain  on state space $S$ with a transition tensor $\mathcal{P}$.  Let $\pi$ be its unique marginal stationary distribution. Then, for every $k \ge m$ and every $i \in S$, 
\be
\label{geometric2}
| \xi_{i}^{(k)} - \pi_i | \le C_1 D^k,
\ee
for any initial distribution ${\cal Z}^{(1)}=[  z^{(1)}_{i_2 \ldots i_m} ]$,  where $C_1=CD^{1-m}$, and $D$ and $C$ are the same constants as appeared in Theorem \ref{geometric1}.  
\end{3}
\bp
Since
$$
\xi_{i}^{(k)} = \sum_{i_2, \ldots, i_m \in S} p_{i i_2 \ldots i_m}^{(k-m+1)}  z^{(1)}_{i_2 \ldots i_m}
$$
and
$$
\sum_{i_2, \ldots, i_m \in S}  z^{(1)}_{i_2 \ldots i_m}  = 1,
$$
we have
\be
\label{xi-minus-pi}
\xi_{i}^{(k)} - \pi_i = \sum_{i_2, \ldots, i_m \in S} (p_{i i_2 \ldots i_m}^{(k-m+1)} - \pi_i) z^{(1)}_{i_2 \ldots i_m}.
\ee
This implies that
\[\begin{array}{rcl}
|\xi_{i}^{(k)} - \pi_i| & \le & \sum_{i_2, \ldots, i_m \in S} |p_{i i_2 \ldots i_m}^{(k-m+1)} - \pi_i| z^{(1)}_{i_2 \ldots i_m} \\
& \le & 
\max_{i_2, \ldots, i_m}  |p_{i i_2 \ldots i_m}^{(k-m+1)} - \pi_i| \sum_{i_2, \ldots, i_m} z^{(1)}_{i_2 \ldots i_m} \\
& = &  \max_{i_2, \ldots, i_m}  |p_{i i_2 \ldots i_m}^{(k-m+1)} - \pi_i | \\
& \le  & CD^{k-m+1} \\
&= & C_1 D^k,
\end{array}
\]
where $C_1 = CD^{1-m}$. 
\ep

For a regular first order Markov chain, the rate of convergence is often studied via the spectral gap between the dominant eigenvalue (which is $1$) and the subdominant eigenvalue of the transition matrix (see, for example, \cite{LevPer17}).  When the reduced first order chain of a regular higher order chain is regular, we can use the spectral gap of the transition matrix $Q$ of the reduced first order chain to obtain the rate of convergence of the original higher order chain. In \cite{HanXu26},  the following result is proved.

\begin{1} 
Suppose that the transition tensor ${\cal P}$ of  an $(m-1)$th order Markov chain $\{X_t\}_{t=1}^{\infty}$ on state space
$S$ is positive, that is, ${\cal P}>0$. Let $Q$ be the transition matrix of its reduced first order chain.  Then the $(m-1)$th power of $Q$ is a positive matrix,  that is,  
$$
Q^{m-1} >0.
$$
This implies that $Q$ itself is regular.
\end{1}

By the Perron-Frobenius theory, if a square matrix $A$ is positive,  then its spectral radius $\rho(A)$ is a simple eigenvalue of $A$ and any other eigenvalue $\lambda$ of $A$ satisfies $|\lambda| < \rho(A)$ (see, for example, \cite{BerPle94}).  In fact, we have the following stronger eigenvalue bound result in this case. 

\begin{1} {(\rm Hopf's inequality, see \cite{Hop63})} Suppose that $A=[a_{ij}]$ is a positive square matrix.  Let 
$$
\gamma = \max _{i,j} a_{ij}, \  \ \ \  \delta = \min_{i,j} a_{ij} .
$$  
If $\lambda$ is an eigenvalue of $A$ other than  $\rho(A)$, then 
$$
|\lambda| \leq \frac{\gamma-\delta}{\gamma+\delta} \rho(A).
$$  
\end{1}

The following theorem gives a bound for the subdominant eigenvalue of $Q$ when ${\cal P}$ is positive.  

\begin{2}
\label{rate3}
Suppose that $\{X_t\}_{t=1}^{\infty}$ is an $(m-1)$th order regular Markov chain  on state space $S$ whose transition tensor $\mathcal{P}$ is positive.  Let $Q=[q_{i_1i_2 \ldots i_{m-1}, j_1 j_2 \ldots j_{m-1}}]$ be the transition matrix of its reduced first order chain. Denote
$Q^k=[q^{(k)}_{i_1i_2 \ldots i_{m-1}, j_1 j_2 \ldots j_{m-1}}]$
 Let 
$$
\alpha = \max _{i_1i_2 \ldots i_{m-1}, j_1 j_2 \ldots j_{m-1}} q^{(m-1)}_{i_1i_2 \ldots i_{m-1}, j_1 j_2 \ldots j_{m-1}}
$$
and
$$
\beta = \min_{i_1i_2 \ldots i_{m-1}, j_1 j_2 \ldots j_{m-1}} q^{(m-1)}_{i_1i_2 \ldots i_{m-1}, j_1 j_2 \ldots j_{m-1}}.
$$  
Then for any eigenvalue of $Q$ other than $\rho(Q)=1$, we have
\be
\label{bd-eig}
|\lambda| \le \left ( \frac{\alpha-\beta}{\alpha+\beta} \right )^{\frac{1}{m-1}}.
\ee 
\end{2}
\bp 
Note that if $\lambda$ is an eigenvalue of $Q$, then $\lambda^{m-1}$ is an eigenvalue of $Q^{m-1}$. Using Hopf's inequality, if $\lambda$ is an eigenvalue of $Q$ other than $\rho(Q)=1$, we have
$$
|\lambda|^{m-1} \leq \frac{\alpha-\beta}{\alpha+\beta},
$$   
Therefore,
$$
|\lambda| \le \left ( \frac{\alpha-\beta}{\alpha+\beta} \right )^{\frac{1}{m-1}}.
$$  
\ep

The bound (\ref{bd-eig}) provides an estimate on the rate of convergence of $\{Q^k \}$. It also gives an estimate on the rate of convergence of ${\cal P}^{(k)}$ by (\ref{PkQk}).

\section{Marginal mixing times}
\label{Mixing-times}
\setcounter{equation}{0}

 For a regular higher order Markov chain  $\{X_t\}_{t=1}^{\infty}$, its reduced first order chain  $\{Y_t\}_{t=1}^{\infty}$ may be reducible.  When $\{Y_t\}_{t=1}^{\infty}$ is reducible, it does not have a unique joint stationary distribution, and therefore, its mixing times are not well defined. However, $\{X_t\}_{t=1}^{\infty}$ has a unique marginal stationary distribution $\pi$ and the distribution of the current state $X_t$ converges to $\pi$. 
 
 In this section, we introduce two types of marginal mixing times, which extend the notions of mixing times for first order Markov chains to higher order chains.  These marginal mixing times capture useful information on  how long it takes for the visible state $X_t$ of a higher order chain  to look stationary, even though the underlying process has memory. They provide a framework for studying the long term behavior of higher order Markov chains when the current state is what matters.\\

We first define marginal mixing times via total variance distance. Recall that the  total variance distance between two probability distributions $\mu=[\mu_1,\ldots, \mu_n]^T$ and $\nu=[\nu_1,\ldots \nu_n]^T$ on $S$  is defined by 
 \be
 \|\mu - \nu\|_{TV} = \frac{1}{2} \sum_{i \in S} |\mu_i - \nu_i|. 
 \ee
If $\nu^{(1)}, \nu^{(2)}, \ldots$ and $\nu$  are probability distributions  on $S$, we say that $\nu^{(k)}$ converges to $\nu$ in total variance as $ k \to \infty$, if 
$$
\lim_{k \to \infty}  \|\nu^{(k)} - \nu\|_{TV} =0.
$$

Define
\be
\label{dk}
d(k) = \max_{i_2, \ldots, i_m \in S} \| p^{(k)}_{\cdot i_2 \ldots i_m} - \pi \|_{TV},
\ee
where $p^{(k)}_{\cdot i_2 \ldots i_m}$ denotes the vector $[p^{(k)}_{1 i_2 \ldots i_m}, \ldots, p^{(k)}_{n i_2 \ldots i_m}]^T$. 

\begin{4}
Suppose that $\{X_t\}_{t=1}^{\infty}$ is an $(m-1)$th order regular Markov chain  on state space $S$ with a transition tensor $\mathcal{P}$.  Let $\pi$ be its unique marginal stationary distribution and $d(k)$ be defined as in  (\ref{dk}). Then there exist constants $C_2>0$ and $0<D<1$ such that 
\be
\label{dk-bound}
d(k) \le C_2D^k.
\ee
Moreover, the sequence $\{d(k)\}$ is decreasing, that is, 
\be
\label{dk-decreasing}
d(k+1) \le d(k).
\ee
\end{4}
\bp
It is easy to see that (\ref{dk-bound}) follows from Theorem \ref{geometric1}. To prove (\ref{dk-decreasing}), suppose that the starting states $(i_{2}^{\ast}, \ldots, i_{m}^{\ast})$ maximizes $d(k+1)$. Then 
\[
\begin{array}{rcl}
d(k+1) & = & \| p^{(k+1)}_{\cdot i_{2}^{\ast} \ldots i_{m}^{\ast}} - \pi \|_{TV} \\
& = & \| \sum_{j \in S} p_{j i_{2}^{\ast} \ldots i_{m}^{\ast}} p^{(k)}_{\cdot j i_{2}^{\ast} \ldots i_{m-1}^{\ast}}    - \pi \|_{TV} \\
& = & \| \sum_{j \in S} p_{j i_{2}^{\ast} \ldots i_{m}^{\ast}}  ( p^{(k)}_{\cdot j i_{2}^{\ast} \ldots i_{m-1}^{\ast}}    - \pi) \|_{TV} \\
& \le &  \max_{i_2, \ldots, i_m \in S} \| p^{(k)}_{\cdot i_2 \ldots i_m} - \pi \|_{TV} \sum_{j \in S} p_{j i_{2}^{\ast} \ldots i_{m}^{\ast}}  \\
&=& d(k).
\end{array}
\]
\ep
\begin{2}
Suppose that $\{X_t\}_{t=1}^{\infty}$ is an $(m-1)$th order regular Markov chain  on state space $S$ with a transition tensor $\mathcal{P}$.  Let $\pi$ be its unique marginal stationary distribution.  Then starting with any initial joint distribution ${\cal Z}^{(1)}=[z^{(1)}_{i_2 \ldots i_m}]$ of $(X_{m-1}, \ldots, X_1)$, the distribution $\xi^{(k)}$ of $X_k$ (for $k \ge m$) satisfies
\be
\| \xi^{(k)}- \pi \|_{TV} \le d(k-m+1).
\ee
\end{2}
\bp
\[
\begin{array}{rcl}
\| \xi^{(k)}- \pi \|_{TV} & =& \displaystyle{  \frac{1}{2} \sum_{i_1 \in S} | \xi^{(k)}_{i_1} - \pi_{i_1}| } \\
& = & \displaystyle{  \frac{1}{2} \sum_{i_1 \in S}  |  \sum_{i_2,\ldots, i_m \in S} p^{(k-m+1)}_{i_1 i_2 \ldots i_m} z^{(1)}_{i_2 \ldots i_m} - \pi_{i_1}| } \\
& = & 
\displaystyle{  \frac{1}{2} \sum_{i_1 \in S}  |  \sum_{i_2,\ldots, i_m \in S} (p^{(k-m+1)}_{i_1 i_2 \ldots i_m} - \pi_{i_1}) z^{(1)}_{i_2 \ldots i_m}| }   \ \ \  \ \  {\rm (by \ (\ref{xi-minus-pi}) } \\
& \le & \displaystyle{  \frac{1}{2} \sum_{i_1 \in S}  \left ( \max_{i_2, \ldots, i_m \in S} |p^{(k-m+1)}_{i_1 i_2 \ldots i_m} - \pi_{i_1}| \sum_{i_2,\ldots, i_m \in S}  z^{(1)}_{i_2 \ldots i_m} \right )} \\
& =& \displaystyle{  \frac{1}{2} \sum_{i_1 \in S}   \max_{i_2, \ldots, i_m \in S} |p^{(k-m+1)}_{i_1 i_2 \ldots i_m} - \pi_{i_1}|  } \\
& = & \displaystyle{    \max_{i_2, \ldots, i_m \in S} \frac{1}{2} \sum_{i_1 \in S}   |p^{(k-m+1)}_{i_1 i_2 \ldots i_m} - \pi_{i_1}|  } \\
& = & d(k-m+1).
\end{array}
\]
\ep
\begin{10} 
\label{TV-mixing}
{\rm
The marginal mixing time via total variance distance is defined as
$$
t_{mix} (\epsilon) = \min \{ t \ge m: d(t) \le \epsilon \}. 
$$ 
}
\end{10}

\ \\

In \cite{Hun06}, Hunter introduced  a notion of expected mixing times for first order Markov chains. We now extend it to higher order chains by defining expected marginal mixing times.  Suppose that  $\{X_t\}_{t=1}^{\infty}$ is an $(m-1)$th order regular Markov chain on state space $S$ with a transition tensor ${\cal P}$. Let $\pi$ be its unique marginal stationary distribution. Since  $\{X_t\}_{t=1}^{\infty}$ is aslo ergodic, it has a unique
mean first passage time $\mu_{i_1 i_2 \ldots i_m}$ for each $ i_1, i_2,  \ldots, i_m \in S$ according to Theorem \ref{MFPT} .   

Let $Y$ be a random variable that takes values on $S$ and has the  distribution $\pi$.   Let $T_{ i_1 i_2 \ldots i_m }$ be the first passage time random variable from states $(i_2,\ldots,i_m)$ to state $i_1$ as defined in (\ref{first-passage}).

Define the marginal mixing time random variable
$$
T= \min \{ j \ge 1 \ : \ X_{  j+m-1} = Y \}.
$$ 

\begin{10} {\rm
\label{EXP-mixing}
 The expected time to marginal mixing starting with $(X_{m-1}, \ldots, X_{1})=(i_{2},\ldots, i_m)$ is defined as
\be
\label{EMMT}
\theta_{i_2 \ldots i_{m}}= E(T \ | \ X_{m-1}=i_2, \ldots, X_{1}=i_{m}).
\ee
}
\end{10}

The following theorem shows that the expected marginal mixing times can be obtained using the mean first passage times. 

\begin{2}
\label{expected-mixing}
\be
\label{theta}
\theta_{i_2 \ldots i_{m}} =  \sum_{i_1 \in S} \mu_{i_1 i_2 \ldots i_m} \pi_{i_1},
\ee
where $\mu_{i_1 i_2 \ldots i_m}$ is the mean first passage time from   $(i_2,\ldots,i_m)$ to state $i_1$. 
\end{2}
\bp
\[
\begin{array}{rcl}
\theta_{i_2 \ldots i_{m}} & = & \displaystyle{E_Y ((T \ | \ X_{m-1}=i_2, \ldots, X_{1}=i_{m}) \ | \ Y) }\\
& = &  \displaystyle{\sum_{i_1 \in S} E ((T \ | \ X_{m-1}=i_2, \ldots, X_{1}=i_{m}) \ | \ Y=i_1) \Pr (Y=i_1)} \\
& =&  \displaystyle{ \sum_{i_1 \in S} E (T_{i_1 i_2 \ldots i_m} \ | \  X_{m-1}=i_2, \ldots, X_{1}=i_{m} ) \pi_{i_1} }\\
& = &  \displaystyle{ \sum_{i_1 \in S} \mu_{i_1 i_2 \ldots i_m} \pi_{i_1}}.
\end{array}
\]
\ep

\begin{5}{\rm
When $m=2$, the expected maginal mixing times are just the expected mixing times. They are a constant, which is a version of the Kemeny constant (See \cite{Hun06}). When $m \ge 3$, the expected marginal mixing times $\theta_{i_2 \ldots i_{m}}$'s are not a constant in general. However, there are special cases in which $\theta_{i_2 \ldots i_{m}}$'s are constant. For example, when $m=3$, if the transition tensor ${\cal P}$ satisfies
$
{\cal P}(:,:,1) =  {\cal P}(:,:,2) = \cdots = {\cal P}(:,:,n),
$
then we have $\theta_{i_2 i_3} \equiv constant$ for all $i_2, i_3 \in S$.  

An interesting question is: Under what conditions, are $\theta_{i_2 \ldots i_{m}}$'s constant? 
}
\end{5}

\end{document}